\documentclass[11pt,reqno,tbtags]{amsart}
\usepackage{mathrsfs}
\usepackage{amssymb}
\usepackage{amsfonts}
\usepackage{stmaryrd}
\pdfoutput=1
\usepackage[colorlinks=true, pdfstartview=FitV, linkcolor=blue,
  citecolor=blue, urlcolor=blue,pagebackref=false]{hyperref}
\usepackage[numeric,initials,nobysame]{amsrefs}
\usepackage{amsmath,amssymb,amsfonts,dsfont}
\usepackage{enumerate,graphicx}
\usepackage{mathrsfs}
\usepackage{appendix}
\numberwithin{equation}{section}

\newtheorem{theorem}{Theorem}[section]
\newtheorem*{theorem*}{Theorem}
\newtheorem{lemma}[theorem]{Lemma}

\newtheorem{proposition}[theorem]{Proposition}

\theoremstyle{definition}{

\newtheorem{definition}{Definition}

}
\theoremstyle{remark}{
\newtheorem{remark}{Remark}
\newtheorem*{remark*}{Remark}

}

\newcommand{\R}{\mathbb R}

\newcommand{\C}{\mathbb C}

\newcommand{\eps}{\varepsilon}
\renewcommand\phi{\varphi}

\newcommand{\vE}{\mathbf{E}}
\newcommand{\vH}{\mathbf{H}}
\newcommand{\vF}{\mathbf{F}}

\newcommand{\deq}{:=}
\newcommand{\Div}{\mbox{Div}}

\newcommand{\Osd}{\Omega\setminus\overline{D}}
\newcommand{\tvE}{\tilde{\vE}}
\newcommand{\vP}{\mathbf{P}}
\newcommand{\vQ}{\mathbf{Q}}
\begin{document}
\title{Reconstructing electromagnetic obstacles by the enclosure method }
\author{Ting Zhou }
\address{Department of Mathematics, University of Washington, Seattle, WA 98195 USA, email: tzhou@math.washington.edu}

\maketitle

\begin{abstract}
We present a reconstruction algorithm for recovering both "magnetic-hard" and "magnetic-soft" obstacles in a background domain
with known isotropic medium from the boundary impedance map. We use in our algorithm complex geometric
optics solutions constructed for Maxwell's equation.
\end{abstract}

\section{Introduction}
In this paper we study an inverse boundary value problem for Maxwell's equation. Let $\Omega\subset\R^3$ be a bounded domain with smooth boundary, filled
with isotropic electromagnetic medium, characterized by three parameters: the permittivity $\eps(x)$, conductivity $\sigma(x)$ and permeability
$\mu(x)$. A "magnetic-hard" obstacle is a subset $D$ of $\Omega$, with smooth boundary, such that the electric-magnetic field $(\vE, \vH)$ satisfies the following BVP for
Maxwell's equation
\begin{equation}
\left\{\begin{array}{l}\nabla\wedge\vE=i\omega\mu\vH,\;\;\;\;\nabla\wedge\vH=-i\omega(\eps+i\frac{\sigma}{\omega})\vE\;\;\;\;\mbox{in }\;\Osd,\\
\nu\wedge\vE|_{\partial\Omega}=f,\\ \nu\wedge\vH|_{\partial D}=0\end{array}\right.
\end{equation}
where $\nu$ is the unit outer normal vector to the boundary $\partial\Omega\cup\partial D$. The boundary condition $\nu\wedge\vH|_{\partial D}=0$ physically appears when an active object is presented. Another type of obstacle, what we call a "magnetic-soft" obstacle, is a subdomain $D$ such that the tangential component of the electric field $\nu\wedge\vE$ vanishes on the interface $\partial D$, when a passive object is presented.  Define the impedance map by taking the tangential component of the
electric field $\nu\wedge\vE|_{\partial\Omega}$ to the tangential component of the magnetic field $\nu\wedge\vH_{\partial \Omega}$. Then our purpose is to
retrieve information of the shape of $D$ from the impedance map.\\

The well-known Calder\'on's problem \cite{C} is to determine the conductivity of a medium by making voltage and current measurements of the boundary. The information is encoded in the Dirichlet-to-Neumann map for the conductivity equation $\nabla\cdot(\gamma\nabla u)=0$. In \cite{SU}, Sylvester and Uhlmann constructed complex geometric optics
(CGO) solutions for Schr\"odinger operator $\Delta-q$ and proved the uniqueness of $C^2$ isotropic conductivity in dimensions $n\geq3$. Further
developments including improved regularity assumption, $2$D problems and partial data problems, were obtained based on this idea. See \cite{U} for a survey of recent developments.

Another application of CGO solutions, the enclosure method was first introduced by Ikehata \cite{I, I2} to identify obstacles, cavities and inclusions embedded in
conductive or acoustic medium. Geometrically, using the property of CGO solutions that decay on one side and grow on the other
side of a hyperplane, one can enclosing obstacles by those hyperplanes. This idea was generalized to identify non-convex obstacles by Ide et al. \cite{IINSU} for isotropic conductivity
equations (conductive medium) and  by Nakamura and Yoshida \cite{NY} for Helmholtz equations (acoustic medium), by utilizing the so-called complex
spherical waves (CSW), namely, CGO solutions with nonlinear Carleman limiting weights. In \cite{UW}, Uhlmann and Wang constructed generalized complex geometric optics solutions for several systems with Laplacian as the leading order term, e.g., the isotropic elasticity system, and implemented them to reconstruct inclusions.  

As for Maxwell's equation,  \cite{SIC}, \cite{OPS1} and \cite{OS} answered the uniqueness question for parameters with suitable regularity from the impedance map in a domain $\Omega$. In \cite{OS}, the Maxwell's operator was reduced into a matrix Schr\"odinger operator and vector CGO solutions were constructed to recover electromagnetic parameters. 

To address the inverse problem of determining an electromagnetic obstacle, we observe that solutions of a non-dissipative Maxwell's equation ($\sigma=0$) share similar asymptotical behavior (a key equality in Lemma \ref{lem_1}) to those of Helmholtz equations (a key inequality in Lemma 4.1 in \cite{I}). Therefore, with CGO solutions at hand, the enclosure method is applicable: one can define an indicator function $I_\rho(\tau, t)$ for each direction $\rho\in \mathbb{S}^2$; by adjusting $t$, the hyperplane moves along $\rho$; for each $\rho$ and $t$, the asymptotical behavior of $I_\rho(\tau, t)$ as $\tau\gg1$ produces the support function of the convex hull of $D$. However, unlike the Schr\"odinger equation, the CGO solution for Maxwell's equation doesn't behave as small perturbation (w.r.t. $\tau$) of Calder\'on's solutions. We overcome this by carefully choosing relatively large "incoming" fields (w.r.t. $\tau$) compared to the perturbation. \\

The rest of the paper is organized as following.  In Section 2, we formulate the forward problem for a "magnetic-hard" obstacle and define rigorously the impedance map. Then we construct the CGO solution of interests in Section 3. The main reconstruction algorithm for "magnetic-hard" obstacles is introduced and proved in Section 4. Finally, we remark in Section 5 that the scheme also applies to "magnetic-soft" obstacles. Through the whole paper, well-posedness of a mixed boundary value problem for Maxwell's equation plays an important role. Hence, for completeness, we include a proof in Appendix A. 

\section{Direct Problems and the main result}
Let $\Omega$ be a bounded domain in $\R^3$ with smooth
boundary and
its complement $\R^3\setminus \bar{\Omega}$ is connected. We consider electric permittivity $\eps(x)$,
conductivity $\sigma(x)$ and
magnetic permeability $\mu(x)$ of the background medium as globally defined functions with following properties:
there are positive constants $\eps_m, \eps_M, \mu_m,
\mu_M$, $\sigma_M$, $\eps_0$ and $\mu_0$ such that for all $x\in \Omega$
\[\eps_m\leq\eps(x)\leq \eps_M,\;\; \mu_m\leq\mu(x)\leq\eps_M,\;\;0\leq\sigma(x)\leq\sigma_M\] and $\eps-\eps_0, \sigma, \mu-\mu_0\in
C_0^3(\Omega)$.\\
A "magnetic-hard" obstacle $D$ (corresponding to the sound-hard obstacle for Helmholtz equations) is a subset of $\Omega$ such that
$\Omega\setminus\overline{D}$ is connected. Moreover, the electric field $\vE$ and the magnetic field $\vH$ in $\Omega\setminus\overline{D}$
satisfy the following boundary value problem of the time-harmonic Maxwell's equation
\begin{equation}\label{Max_1}
\left\{\begin{array}{l}\nabla\wedge \vE=i\omega\mu \vH,\;\;\nabla \wedge \vH=-i\omega\gamma\vE \;\;\;\;\mbox{in}\;\;\Osd,\\
\nu\wedge\vE=f\in TH^{1/2}_{\mbox{Div}}(\partial\Omega) \;\;\;\;\mbox{on}\;\;\partial\Omega,
\end{array}\right.
\end{equation}
where $\gamma=\eps+i\frac{\sigma}{\omega}$,
and the "magnetic-hard" boundary condition on the interface $\partial D$
\begin{equation}\label{bd_1}(\nu\wedge\vH)|_{\partial D}=0,\end{equation} where $\nu$ is the unit outer normal vector on $\partial\Omega\cup\partial D$.
Through out this note, we assume the non-dissipative case $\sigma=0$. Then $\gamma=\eps$ is a real function. \\

\noindent \textbf{Notations.} If $F$ is a function space on $\partial \Omega$, the subspace of all those $f\in F^3$ which are tangent to
$\partial \Omega$ is denoted by $TF$. For example,
for $u\in (H^s(\partial\Omega))^3$, we have decomposition $u=u_t+u_\nu\nu$, where the tangential component
$u_t=-\nu\wedge(\nu\wedge u)\in TH^s(\partial\Omega)$ and the normal component $u_\nu=u\cdot \nu\in H^s(\partial\Omega)$.
For a bounded domain $\Omega$ in $\R^3$, we denote
\[TH^s_{\Div}(\partial\Omega)=\{f\in TH^s(\partial\Omega)\;|\;\mbox{Div}(f)\in H^{s}(\partial\Omega)\},\]
\[H^k_{\Div}(\Omega)=\{u\in (H^k(\Omega))^3\;|\;\mbox{Div}(\nu\wedge u|_{\partial\Omega})\in H^{k-1/2}(\partial\Omega)\},\] with norms
\[\|f\|^2_{TH^s_{\Div}(\partial\Omega)}=\|f\|^2_{H^s(\partial\Omega)}+\|\mbox{Div}(f)\|^2_{H^{s}(\partial\Omega)},\]
\[\|u\|^2_{H^k_{\Div}(\partial\Omega)}=\|u\|^2_{H^k(\Omega)}+\|\mbox{Div}(\nu\wedge u|_{\partial\Omega})\|^2_{H^{k-1/2}(\partial\Omega)},\]
where Div is the surface divergence. There are natural inner products making them Hilbert spaces (see \cite{SIC}).
In addition, we define the weighted $L^2$ space in $\R^3$:
\[L^2_{\delta}=\left\{f\in L^2_{loc}(\R^3)\;:\;\|f\|_{L^2_{\delta}}^2=\int (1+|x|^2)^\delta|f(x)|^2dx<\infty \right\}.\]\\

\noindent\textbf{Admissibility.} It can be shown (see Appendix A.) that for $f\in TH_{\Div}^{1/2}(\partial\Omega)$ and $g\in TH_{\Div}^{1/2}(\partial D)$, the boundary value problem for
Maxwell's equation
\begin{equation}\label{Max_4}
\left\{
\begin{array}{l}
\nabla\wedge \vE=i\omega\mu \vH,\;\;\nabla \wedge \vH=-i\omega\gamma\vE \;\;\;\mbox{in}\;\;\Omega\setminus\bar{D},\\
\nu\wedge\vE|_{\partial\Omega}=f\\
\nu\wedge\vH|_{\partial D}=g,
\end{array}
\right.
\end{equation}
has a unique solution $(\vE, \vH)\in H^1_{\Div}(\Osd)\times H^1_{\Div}(\Osd)$,
 except for a discrete set of magnetic resonance frequencies $\{\omega_n\}$.
It satisfies
\begin{equation}\|\vE\|_{H^1_{\Div}(\Osd)}+\|\vH\|_{H^1_{\Div}(\Osd)}\leq C\left(\|f\|_{TH^{1/2}_{\Div}(\partial \Omega)}+\|g\|_{TH^{1/2}_{\Div}(\partial D)}\right).\end{equation}

Denote by $(\vE_0, \vH_0)$ the solution of Maxwell's equation in the domain $\Omega$ without an obstacle, namely, the solution of
\[\left\{\begin{array}{l}
\nabla\wedge\vE_0=i\omega\mu\vH_0,\;\;\;\nabla\wedge\vH_0=-i\omega\eps\vE_0,\;\;\;\mbox{in}\;\Omega,\\
\nu\wedge\vE_0|_{\partial\Omega}=f.\end{array}
\right.\]\\

\noindent\textbf{Main result.} Now we are in the position to define the impedance map for non-resonant frequencies, 
\[\Lambda_D(\nu\wedge \vE|_{\partial\Omega})=\nu\wedge \vH|_{\partial\Omega},\] and
it can be shown that \[\Lambda_D:\;TH^{1/2}_{\Div}(\partial\Omega)\rightarrow TH^{1/2}_{\Div}(\partial\Omega)\] is bounded. If $\omega$ is a resonance frequency, one can replace the impedance map by the Cauchy data set 
\begin{eqnarray*}\mathcal{C}_\omega&=&\left\{\left(\nu\wedge\vE|_{\partial\Omega}, \nu\wedge\vE|_{\partial\Omega}\right)\,|\, (\vE, \vH) \mbox{ satisfies } \eqref{Max_1} \mbox{ and } \eqref{bd_1}\right\}\\ &\subset& TH^{1/2}_{\Div}(\partial\Omega)\times TH^{1/2}_{\Div}(\partial\Omega).\end{eqnarray*}
Denote by $\Lambda_\emptyset$ the impedance map for the domain without an obstacle. Then the main result of the presenting work is to show
\begin{theorem}\label{thm_1}
For non-dissipative background medium ($\sigma=0$), there exists a reconstruction scheme for the obstacle $D$ from the impedance map $\Lambda_D$.
\end{theorem}

\section{Complex geometric optics solutions} In \cite{OS}, the Maxwell's operator was reduced to an $8\times 8$ second order Schr\"{o}dinger matrix operator by
introducing the generalized Sommerfeld potential. A vector CGO-solution was constructed for the Schr\"{o}dinger operator,
simplifying the proof in \cite{OPS1}. Similar techniques also appeared in \cite{COS} when dealing with the inverse boundary value
problems for Maxwell's equations with partial data.
For completeness, we include the construction of the solution in this work (see \cite{OS} for more details) and provides a special choice of the "incoming" fields.\\

Define the scalar fields $\Phi$ and $\Psi$ as
\begin{equation}\label{scf}\Phi=\frac{i}{\omega}\nabla\cdot(\gamma\vE),\;\;\;\Psi=\frac{i}{\omega}\nabla\cdot(\mu\vH).\end{equation} Under appropriate assumptions on $\Phi$ and $\Psi$,
Maxwell's equation is equivalent to
\begin{equation}\label{scf_eqn}\nabla\wedge\vE-\frac{1}{\gamma}\nabla\left(\frac{1}{\mu}\Psi\right)-i\omega\mu\vH=0,\;\;\nabla\wedge\vH+\frac{1}{\mu}\nabla
\left(\frac{1}{\gamma}\Phi\right)+i\omega\gamma\vE=0.\end{equation} Moreover, in this case, $\Phi$ and $\Psi$ vanish, leading to a solution of Maxwell's equation. Let $X=(\phi, e, h, \psi)^T\in
(\mathcal{D}')^8$ with
\[e=\gamma^{1/2}\vE, \;\;\;h=\mu^{1/2}\vH,\]
\[\phi=\frac{1}{\gamma\mu^{1/2}}\Phi,\;\;\; \psi=\frac{1}{\gamma^{1/2}\mu}\Psi.\]
Then \eqref{scf} and \eqref{scf_eqn} read
\begin{equation}(P(i\nabla)-k+V)X=0,\;\;\;\mbox{in }\Omega
\end{equation}
where \[P(i\nabla)=\left(\begin{array}{cccc}0 & \nabla\cdot & 0 & 0 \\\nabla & 0 & \nabla\wedge & 0 \\0 & -\nabla\wedge & 0 & \nabla \\
0 & 0 & \nabla\cdot & 0\end{array}\right),\]
\[V=(k-\kappa)\mathbf{1}_8+\left(\left(\begin{array}{cccc}0 & \nabla\cdot & 0 & 0 \\\nabla & 0 & -\nabla\wedge & 0 \\
0 & \nabla\wedge & 0 & \nabla \\0 & 0 & \nabla\cdot & 0\end{array}\right)D\right)D^{-1}\] are matrix operators and
\[D=\mbox{diag}(\mu^{1/2}, \gamma^{1/2}\mathbf{1}_3, \mu^{1/2}\mathbf{1}_3, \gamma^{1/2}),\;\;
\kappa=\omega(\gamma\mu)^{1/2},\;\;k=\omega(\eps_0\mu_0)^{1/2}.\] An important property of this operator is
that it allows to reduce Maxwell's equation to the Schr\"odinger matrix equation by noticing
\begin{equation}(P(i\nabla)-k+V)(P(i\nabla)+k-V^T)=-(\Delta+k^2)\mathbf{1}_8+Q,\end{equation} where
\[Q=VP(i\nabla)-P(i\nabla)V^T+k(V+V^T)-VV^T\] is a zeroth-order matrix multiplier. Hence, by writing an ansatz for $X$,
one can define the generalized Sommerfeld potential $Y$
\[X=(P(i\nabla)+k-V^T)Y.\]
So it satisfies the Schr\"{o}dinger equation \begin{equation}\label{SP}(-\Delta-k^2+Q)Y=0.\end{equation} The following CGO-solution is constructed using Faddeev's Kernel. Let $\zeta\in \C^3$ be a vector with $\zeta\cdot\zeta=k^2$. Suppose $y_{0, \zeta}\in \C^8$ is a constant vector with
respect to $x$ and bounded with respect to $\zeta$. We refer $e^{ix\cdot\zeta}y_{0,\zeta}$ as the "incoming" field. Then there exists a unique solution of  \eqref{SP} of the form
\[Y_\zeta(x)=e^{ix\cdot\zeta}(y_{0, \zeta}-v_\zeta(x)),\] where $v_\zeta(x)\in (L^2_{\delta+1})^8$
satisfying \[\|v_\zeta\|_{L^2_{\delta+1}}\leq C/|\zeta|\] for $\delta\in (-1, 0)$. Moreover, one can show
that $v_\zeta\in ({H^s(\Omega)})^8$ for $0\leq s\leq2$, e.g., see \cite{BU}, and
\begin{equation}\label{remainder_1}
\|v_\zeta(x)\|_{H^s(\Omega)}\leq C|\zeta|^{s-1}.
\end{equation}\\

Lemma 3.1 in \cite{OS} states that if we choose $y_{0, \zeta}$ such that the first and the last components of $(P(\zeta)-k)y_{0, \zeta}$ vanish, where $P(\zeta)$ is the matrix obtained by replacing $i\nabla$ by $\zeta$ in $P(i\nabla)$,
then for large $|\zeta|$, $X_\zeta$ provides the solution of the original Maxwell's equation.
We proceed to provide a more specific choice of $y_{0,\zeta}$ such that the CGO solution for Maxwell's equation has special properties.  \\

As in \cite{OS}, choose
\[y_{0, \zeta}=\frac{1}{|\zeta|}(\zeta\cdot a, ka, kb, \zeta\cdot b)^T,\] where
\[\zeta=-i\tau\rho+\sqrt{\tau^2+k^2}\rho^\bot,\] with $\rho, \rho^\bot\in \mathbb{S}^2$ and $\rho\cdot\rho^\bot=0$. $\tau>0$ is used to control the size of $|\zeta|=\sqrt{2\tau^2+k^2}$.
Then we obtain
\[x_{0,\zeta}:=(P(-\zeta)+k)y_{0, \zeta}=\frac{1}{|\zeta|}\left(\begin{array}{c}0 \\
-(\zeta\cdot a)\zeta-k\zeta\wedge b+k^2 a \\k\zeta\wedge a-(\zeta\cdot b)\zeta+k^2 b \\0\end{array}\right)\] satisfying the condition in Lemma 3.1 in \cite{OS}.\\
Taking $\tau\rightarrow\infty$, we have
\[\frac{\zeta}{|\zeta|}\rightarrow \hat{\zeta}=\frac{1}{\sqrt{2}}(-i\rho+\rho^\bot).\]
We choose $a$ and $b$ such that
\[\hat{\zeta}\cdot b=1,\;\;\;\hat{\zeta}\cdot a=0.\] This is satisfied, for example, by taking $a\in \R^3$ and $b\in \C^3$ satisfying
$a\perp \rho$, $a\perp \rho^\bot$ and $b=\overline{\hat{\zeta}}$. Given these choices,
It's easy to see that \[\eta\deq(x_{0,\zeta})_2\rightarrow -k\hat{\zeta}\wedge b=ik\rho\wedge\rho^\bot\;\; (\sim \mathcal{O}(1))\;\;\mbox{as }\tau\rightarrow\infty,\]
\[\theta\deq(x_{0,\zeta})_3\sim \mathcal{O}(\tau)\;\;\mbox{as }\tau\rightarrow\infty.\] Then $X_\zeta=(P(i\nabla)+k-V^T)Y_\zeta$ is written in the form
\[X_{\zeta}=e^{\tau(x\cdot\rho)+i\sqrt{\tau^2+k^2}x\cdot\rho^\bot}(x_{0,\zeta}+r_\zeta(x))\]
where \begin{equation}\label{remainder_2}r_\zeta=P(-\zeta)v_\zeta+P(i\nabla)v_\zeta-V^Ty_{0,\zeta}+kv_\zeta-V^Tv_\zeta\end{equation} satisfying for $C>0$ independent of $\zeta$
\[\|r_\zeta\|_{L^2(\Omega)}\leq C.\]
Summing up, we obtain the following
\begin{proposition}\label{prop_1}
Let $\rho, \rho^\bot\in \mathbb{S}^2$ with $\rho\cdot\rho^\bot=0$. Assume $\omega$ is not a resonant frequency. Given $\theta, \eta\in \C^3$ as above, then for $\tau>0$ large enough, there exists a unique complex geometric optics solution
$(\vE_0, \vH_0)$ of Maxwell's equation
\[\nabla\wedge\vE_0=i\omega\mu\vH_0\;\;\;\;\nabla\wedge\vH_0=-i\omega\eps\vE_0\;\;\;\;\mbox{in }\;\Omega\] of the form
\[\begin{array}{c}\vE_0=\eps^{-1/2}e^{\tau(x\cdot\rho)+i\sqrt{\tau^2+k^2}x\cdot\rho^\bot}(\eta+R)\\
\vH_0=\mu^{-1/2}e^{
\tau(x\cdot\rho)+i\sqrt{\tau^2+k^2}x\cdot\rho^\bot}(\theta+Q).\end{array}\]
Moreover, we have \[\eta=\mathcal{O}(1),\;\; \theta=\mathcal{O}(\tau)\;\;\;\mbox{for}\;\tau\gg1,\] and
$R=(r_\zeta)_2, Q=(r_\zeta)_3$ are bounded in $(L^2(\Omega))^3$ for $\tau\gg1$.
\end{proposition}

For reconstruction, one needs to compute the boundary tangential CGO-fields $(\nu\wedge\vE_0|_{\partial\Omega}, \nu\wedge\vH_0|_{\partial\Omega})$. In \cite{OS, OPS1}, by solving a boundary integral equation, one can recover them from the impedance map $\Lambda_\emptyset$ if the background parameters are unknown. In our case with known medium, CGO-fields are known.

\section{Reconstruction Scheme}
Adding a parameter $t>0$ into the CGO-solution in Proposition \ref{prop_1}, we use
\begin{equation}\label{CGO}\begin{array}{c}\vE_0=\eps^{-1/2}e^{\tau(x\cdot\rho-t)+i\sqrt{\tau^2+k^2}x\cdot\rho^\bot}(\eta+R)\\ \vH_0=\mu^{-1/2}e^{
\tau(x\cdot\rho-t)+i\sqrt{\tau^2+k^2}x\cdot\rho^\bot}(\theta+Q)\end{array}\end{equation} to define an indicator function which physically measures the differences between the energies required to keep the same boundary CGO electric field for the domain $\Omega$ with and without the obstacle $D$.
\begin{definition}For $\rho\in \mathbb{S}^2$, $\tau>0$ and $t>0$ we define the indicator function \[I_{\rho}(\tau, t)\deq\int_{\partial\Omega}(\nu\wedge\vE_0)\cdot\left(\overline{(\Lambda_D-\Lambda_\emptyset)
(\nu\wedge\vE_0)\wedge\nu}\right)dS\] where $\vE_0$ is a CGO solution of Maxwell's equation given by \eqref{CGO}.
\end{definition}
The enclosure method's aim is to recover the convex hull $\overline{\mbox{ch}(D)}$ of $D$ by reconstructing the following support function.
\begin{definition} For $\rho\in\mathbb{S}^2$, we define the support function of $D$ by
\[h_D(\rho):=\sup_{x\in D}x\cdot\rho.\]
\end{definition}

Then, the reconstruction scheme in Theorem \ref{thm_1} is
\begin{theorem}\label{thm_2}
We assume that the set $\{x\in\R^3\;|\;x\cdot\rho=h_D(\rho)\}\cap\partial D$ consists of one point and the Gaussian
curvature of $\partial D$ is not vanishing at that point. Then, we can recover $h_D(\rho)$ by
\[h_D(\rho)=\inf \{t\in\R\;|\;\lim_{\tau\rightarrow\infty}I_\rho(\tau, t)=0\}.\]
Moreover, if $D$ is strictly convex (the Gaussian curvature is everywhere positive), then we can reconstruct $D$.
\end{theorem}
\begin{remark}\label{rk_1}
The proof of Theorem \ref{thm_2} mainly consists of showing the following limits:
\[
\lim_{\tau\rightarrow\infty}I_\rho(\tau, t)=0,\;\;\;\mbox{when}\;\;t>h_D(\rho);
\]
\[
\liminf_{\tau\rightarrow\infty}I_\rho(\tau, h_D(\rho))=C >0.
\]
\end{remark}
\begin{remark}\emph{Gaussian curvature} at a point on a surface is defined to be the product of two principal curvatures, which measure how the surface bends by different amounts in different directions at that point. A surface with positive Gaussian curvature at a point is locally convex. Note that the non-vanishing assumption on the Gaussian curvature in the theorem is not crucial since only finitely many directions $\rho$ violate the condition.  
\end{remark}

\subsection{A key integral equality}
To show the limits in Remark \ref{rk_1}, we need the following equality for \emph{non-dissipative} Maxwell's equation
\begin{lemma}\label{lem_1}Let $\sigma=0$.
Assume $(\vE, \vH)$ is a solution of
\[\nabla\wedge\vE=i\omega\mu\vH,\;\;\;\;\nabla\wedge\vH=-i\omega\eps\vE,\;\;\;\;\mbox{in }\;\Osd\] satisfying the boundary condition
\[\nu\wedge\vH|_{\partial D}=0\;\;\;\mbox{and}\;\;\;\nu\wedge\vE|_{\partial\Omega}=\nu\wedge\vE_0|_{\partial\Omega}.\] Then we have
\begin{align}\label{Iden_1}
&i\omega\int_{\partial\Omega}(\nu\wedge\vE_0)\cdot\left[\overline{(\nu\wedge\vH-\nu\wedge\vH_0)}\wedge\nu\right]dS\nonumber\\
&=\int_{\Omega\setminus\bar{D}}\mu^{-1}|\nabla\wedge \vE-\nabla\wedge\vE_0|^2-\omega^2\eps|\vE-\vE_0|^2dx\nonumber\\
&+\int_D\mu^{-1}|\nabla\wedge\vE_0|^2-\omega^2\eps|\vE_0|^2dx.
\end{align}
\end{lemma}
\textbf{Proof:} Denote
\[I:=i\omega\int_{\partial\Omega}(\nu\wedge \vE_0)\cdot[\overline{(\nu\wedge \vH-\nu\wedge\vH_0)\wedge\nu }]dS=
i\omega\int_{\partial\Omega}(\nu\wedge\vE_0)\cdot(\overline{\vH-\vH_0})dS.\]
First by integration by parts, we have
\begin{eqnarray*}
&&\int_{\Omega\setminus\bar{D}}\mu^{-1}(\nabla\wedge\vE)\cdot(\overline{\nabla\wedge\vE-\nabla\wedge\vE_0})
-\omega^2\eps\vE\cdot(\overline{\vE-\vE_0})dx\\
&=&-\left(\int_{\partial\Omega}-\int_{\partial D}\right)(\nu\wedge\mu^{-1}(\nabla\wedge\vE))\cdot(\overline{\vE-\vE_0})dS=0
\end{eqnarray*}
by boundary conditions.
Adding this to the following equality
\begin{eqnarray*}
I&=&\int_{\partial\Omega}(\nu\wedge\vE_0)\cdot(\overline{-i\omega\vH+i\omega\vH_0})dS\\
&=&\int_{\Omega\setminus\bar{D}}-\mu^{-1}(\nabla\wedge\vE_0)\cdot(\overline{\nabla\wedge\vE})+\omega^2\eps\vE_0\cdot\overline{\vE}dx\\
&&+\int_\Omega\mu^{-1}|\nabla\wedge\vE_0|^2-\omega^2\eps|\vE_0|^2dx+\int_{\partial D}(\nu\wedge\vE_0)\cdot(\overline{-i\omega\vH})dS
\end{eqnarray*}
with the last term vanishing due to the zero-boundary condition on the interface,
\[
\int_{\partial D}(\nu\wedge\vE_0)\cdot(\overline{-i\omega\vH})dS=\int_{\partial D}(\nu\wedge\vE_0)
\cdot(\overline{-i\omega(\nu\wedge\vH)\wedge\nu})dS=0,
\]
we obtain \eqref{Iden_1}. \qed\\

\subsection{Proof of Theorem \ref{thm_2}}
We proceed to show the first limit in Remark \ref{rk_1}
\begin{equation}\label{lim_1}
\lim_{\tau\rightarrow\infty}I_\rho(\tau, t)=0\;\;\;\mbox{if }\; t>h_D(\rho)
\end{equation}
 by proposing an upper bound of the indicator function.\\

Let $\tilde{\vE}=\vE-\vE_0$ be the reflected solution in $\Osd$. It satisfies
\begin{equation}
\left\{
\begin{array}{l}
\nabla\wedge(\mu^{-1}\nabla\wedge \tvE)-s\eps^{-1}\nabla(\nabla\cdot\eps\tvE)-\omega^2\eps\tvE=0\;\;\;\mbox{in }\;\Osd,\\
\nu\wedge\tvE|_{\partial\Omega}=0,\\
\nu\wedge(\mu^{-1}\nabla\wedge\tvE)|_{\partial D}=-i\omega\nu\wedge \vH_0|_{\partial D}\in TH_{\Div}^{1/2}(\partial D).
\end{array}
\right.
\end{equation}
By \eqref{est_A1} in Appendix A., we have
\begin{equation}\label{tvE_1}\|\tvE\|^2_{H^1(\Osd)}\leq C\|\nu\wedge\vH_0\|^2_{H^{-1/2}(\partial D)}\leq C(\|\nabla\wedge\vE_0\|^2_{L^2(D)}+\|\vE_0\|^2_{L^2(D)}),\end{equation} where the second inequality is valid since
\[i\omega\langle \nu\wedge\vH_0, \vF\rangle_{\partial D}=(\mu^{-1}\nabla\wedge\vE_0, \nabla\wedge\vF)_D-(i\omega\eps\vE_0, \vF)_D\;\;\;\mbox{for}\;F\in (H^1(D))^3.\]
Therefore, by \eqref{Iden_1} and \eqref{tvE_1}, we have
\begin{equation}
I_\rho(\tau, t)\leq C(\|\vE_0\|^2_{L^2(D)}+\|\nabla\wedge\vE_0\|^2_{L^2(D)})\leq C(\|\vE_0\|^2_{L^2(D)}+\|\vH_0\|^2_{L^2(D)}).
\end{equation}
Plugging in the CGO-solution \eqref{CGO}, we obtain the following estimates:
\begin{equation}\label{est_E0_H0}\begin{array}{c}\|\vE_0\|^2_{{L^2(D)}}\leq Ce^{2\tau(h_D(\rho)-t)}\|\eta+R\|^2_{{L^2(D)}^3}\sim e^{2\tau(h_D(\rho)-t)}\;\;\;\tau\gg1,\\
\|\vH_0\|^2_{{L^2(D)}}\leq Ce^{2\tau(h_D(\rho)-t)}\|\theta+Q\|^2_{{L^2(D)}^3}\sim \tau^2e^{2\tau(h_D(\rho)-t)}\;\;\;\tau\gg1.\end{array}\end{equation}
Therefore, we obtain
\[I_\rho(\tau, t)\leq C\tau^2  e^{2\tau(h_D(\rho)-t)}\] for $\tau$ large enough, proving the first limit \eqref{lim_1}.\\

To show the second limit
\begin{equation}\label{lim_2}\liminf_{\tau\rightarrow\infty}I_\rho(\tau, h_D(\rho))=C >0,\end{equation}
it suffices to prove the following two lemmas.
\begin{lemma}\label{lem_2}
If $t=h_D(\rho)$ in CGO-solution \eqref{CGO}, then
\[\liminf_{\tau\rightarrow\infty}\int_D\mu^{-1}|\nabla\wedge\vE_0|^2dx=C,\]with some constant $C>0$.
\end{lemma}
\begin{lemma}\label{lem_3}
If $t=h_D(\rho)$, then there exists a positive number $c$ such that
\begin{equation}\label{lem3_1}\frac{\omega^2\left(\int_{\Osd}\eps|\vE-\vE_0|^2dx+\int_D\eps|\vE_0|^2dx\right)}{\int_D\mu^{-1}|\nabla\wedge\vE_0|^2dx}\leq c<1,\end{equation}
for $\tau$ large enough.
\end{lemma}

\textbf{Proof of Lemma \ref{lem_2}}: This is obtained by noticing, in Proposition \ref{prop_1}, that the first order growth of the constant vector
$\theta$ in $\vH_0$ with respect to $\tau$.
Then the left hand side integral
\[\int_D\mu^{-1}|\nabla\wedge\vE_0|^2dx\geq C\|\vH_0\|^2_{{L^2(D)}^3}\geq C\int_D\tau^2e^{2\tau(x\cdot\rho-h_D(\rho))}dx\geq C\;\;\mbox{for} \;\tau\gg1.\]
To show the last inequality, we denote by $x_0$ the point in $\{x\in \R^3\,|\,x\cdot\rho=h_D(\rho)\}\cap\partial D$. It's not hard to see that there exist $C_\rho>0$ and $\delta_\rho>0$ such that
\[\mu_2\left(D_\rho(\delta_\rho, s)\right)\geq C_\rho s\] where $\mu_2$ denotes the two dimensional Lebesgue measure and
\[D_\rho(\delta_\rho, s)=\{x\in D\,|\,x\cdot\rho=h_D(\rho)-s\}.\]
Then we decompose $D$ into
\[D_\rho(\delta_\rho)=\{x\in D\,|\,h_D(\rho)-\delta_\rho<x\cdot\rho\leq h_D(\rho)\}\] and $D\setminus D_\rho(\delta_\rho)$. The integral of $\tau^2e^{2\tau(x\cdot\rho-h_D(\rho))}$ on $D\setminus D_\rho(\delta_\rho)$ vanishes as $\tau\rightarrow\infty$.  On $D_\rho(\delta_\rho)$, we have
\begin{eqnarray*}
\int_{D_\rho(\delta_\rho)}\tau^2e^{2\tau(x\cdot\rho-h_D(\rho))}&=&\tau^2\int_0^{\delta_\rho}ds\int_{D_\rho(\delta_\rho, s)}e^{-2\tau s}dS\\
&\geq&C_\rho\tau^2\int_0^{\delta_\rho}se^{-2\tau s}ds\\
&=&C_\rho\int_0^{\tau\delta_\rho}se^{-2s}ds\rightarrow \frac{1}{4}C_\rho\;\;\;\;\;\mbox{as}\;\tau\rightarrow\infty.
\end{eqnarray*}\qed\\

\textbf{Proof of Lemma \ref{lem_3}}: The proof follows a similar scheme to \cite{I} for Helmholtz equations.\\
Noticing that
\[\frac{\omega^2\int_D\eps|\vE_0|^2dx}{\int_D\mu^{-1}|\nabla\wedge\vE_0|^2dx}=\frac{\int_D\eps|\vE_0|^2dx}{\int_D\mu|\vH_0|^2dx} \sim \mathcal{O}(\tau^{-2})\;\;\;\;\mbox{for }\;\tau\gg1,\]
it suffices to show
\[\lim_{\tau\rightarrow\infty}\frac{\omega^2\int_{\Osd}\eps|\vE-\vE_0|^2dx}{\int_D\mu^{-1}|\nabla\wedge\vE_0|^2dx}
=\lim_{\tau\rightarrow\infty}\frac{\int_{\Osd}\eps|\vE-\vE_0|^2dx}{\int_D\mu|\vH_0|^2dx}=0.\]
To estimate the numerator, we consider the boundary value problem:
\begin{equation}
\left\{\begin{array}{l}
\nabla\wedge \vP=i\omega\mu \vQ,\;\;\nabla\wedge \vQ=-i\omega\eps \vP-\frac{i\eps}{\omega}\overline{\tilde{\vE}}\;\;\;\mbox{in }\;\Osd,\\
\nu\wedge \vP|_{\partial\Omega}=0,\\
\nu \wedge \vQ|_{\partial D}=0.
\end{array}\right.
\end{equation}
or equivalently
\begin{equation}\label{auxbvp}
\left\{\begin{array}{l}
\nabla\wedge(\mu^{-1}\nabla\wedge \vP)-s\eps^{-1}\nabla(\nabla\cdot\eps \vP)-\omega^2\eps\vP=\eps\overline{\tvE}\;\;\mbox{ in }\Osd,\\
\nu\wedge \vP|_{\partial\Omega}=0,\\
\nu\wedge(\nabla\wedge\vP)|_{\partial D}=0.
\end{array}\right.
\end{equation}
Note that \[\nabla\cdot\eps\vP=0,\;\;\;\;\mbox{ in }\Osd\] because $\nabla\cdot\eps\tvE=0$ in $\Osd$.

Since $\omega$ is admissible, the boundary value problem \eqref{auxbvp}
is well-posed for $\eps\overline{\tvE}=\eps(\overline{\vE-\vE_0})$ which is in $H^1_{\Div}(\Osd)$. Moreover, by Proposition \ref{prop_3}, one has $\vP\in (H^2(\Osd))^3$ satisfying
\[\|\vP\|_{H^2(\Osd)}\leq C\|\tvE\|_{L^2(\Osd)}.\]
By the Sobolev embedding theorem, we have
\[|\vP(x)-\vP(y)|\leq C|x-y|^{1/2}\|\tilde{\vE}\|_{L^2(\Osd)}\;\;\;\mbox{for }\;x, y\in \Osd,\]
\[\sup_{x\in\Osd}|\vP(x)|\leq C\|\tilde{\vE}\|_{L^2(\Osd)}.\]
Since \[\nabla\wedge(\mu^{-1}\nabla\wedge \vP)-\omega^2\eps \vP=\eps\overline{\tilde{\vE}},\] integration by parts gives
\begin{eqnarray*}
\int_{\Osd}\eps|\tilde{\vE}|^2dx&=&\int_{\Osd}\tilde{\vE}\cdot(\nabla\wedge(\mu^{-1}\nabla\wedge \vP)-\omega^2\eps \vP)dx\\
&=&\int_{\Osd}\mu^{-1}(\nabla\wedge \tilde{\vE})\cdot(\nabla\wedge \vP)-\omega^2\eps\tilde{\vE}\cdot \vP dx\\
&&+\left(\int_{\partial\Omega}-\int_{\partial D}\right)\tilde{\vE}\cdot(\nu\wedge(\mu^{-1}\nabla\wedge \vP))dS\\
&=&\int_{\Osd}\nabla\wedge(\mu^{-1}\nabla\wedge\tilde{\vE})\cdot \vP-\omega^2\eps\tilde{\vE}\cdot \vP dx\\
&&-\left(\int_{\partial\Omega}-\int_{\partial D}\right)\nu\wedge(\mu^{-1}\nabla\wedge\tilde{\vE})\cdot \vP dS\\
&=&-\int_{\partial D}\nu\wedge(\mu^{-1}\nabla\wedge\vE_0)\cdot \vP dS.
\end{eqnarray*}
Expanding the RHS of the last equality at $x_0$, one has
\begin{eqnarray*}
\int_{\Osd}\eps|\tvE|^2dx&=&\int_{\partial D}(\vP(x_0)-\vP(x))\cdot \nu\wedge(\mu^{-1}\nabla\wedge\vE_0)dS-\int_D\omega^2\eps\vE_0\cdot\vP(x_0)dx\\
&\leq&C\left\{\int_{\partial D}|x-x_0|^{1/2}|\nu\wedge\vH_0|dS+\int_D|\vE_0|dx\right\}\|\tilde{\vE}\|_{L^2(\Osd)}\\
&\leq&C \left\{\int_{\partial D}\tau|x-x_0|^{1/2}e^{\tau(x\cdot\rho-h_D(\rho))}dS+\int_De^{\tau(x\cdot\rho-h_D(\rho))}dx\right\}\|\tilde{\vE}\|_{L^2(\Osd)}\\
\end{eqnarray*}
This yields
\begin{eqnarray*}
\|\tilde{\vE}\|^2_{L^2(\Osd)}dx
&\leq&C\left(\frac{\int_{\Osd}\eps|\tvE|^2dx}{\|\tilde{\vE}\|_{L^2(\Osd)}}\right)^2\\
&\leq& C\left\{\tau^2\left(\int_{\partial D}|x-x_0|^{1/2} e^{\tau(x\cdot\rho-h_D(\rho))}dS\right)^2+\left(\int_De^{\tau(x\cdot\rho-h_D(\rho))}dx\right)^2\right\}.\end{eqnarray*}
It's easy to see that the second term $\left(\int_D e^{\tau(x\cdot\rho-h_D(\rho))}dx\right)^2$ can be absorbed by the denominator $\int_D\mu^{-1}|\nabla\wedge\vE_0|^2dx$  in \eqref{lem3_1} for large $\tau$. Therefore, it's sufficient to show
\[\lim_{\tau\rightarrow \infty}\tau\int_{\partial D}|x-x_0|^{1/2} e^{\tau(x\cdot\rho-h_D(\rho))}dS=0.\]
This is shown in \cite{I}, where the assumption that the Gaussian curvature of $\partial D$ at $x_0$ is non-vanishing was used. 
This completes the proof of the lemma, hence proves the theorem.\qed

\section{Enclosing "magnetic-soft" obstacles and inclusions}
In  \cite{I2}, the reconstruction procedure for sound-hard obstacles also works for sound-soft obstacles. Inspired by this, our method also applies to enclosing a "magnetic-soft" obstacle.
Suppose our domain $\Omega$, obstacle $D$ and all the electromagnetic parameters in the background satisfy the same hypothesis in the "magnetic-hard" case, except that the fields $(\vE, \vH)$ satisfy Maxwell's equation \eqref{Max_1} with the boundary condition
\begin{equation}\label{soft}\nu\wedge\vE|_{\partial D}=0.\end{equation} Then the reconstruction scheme Theorem \ref{thm_2} applies simply by noticing the following key equality.
\begin{lemma}\label{lem_soft}
Let $\sigma=0$.
Assuming Maxwell's equation with the boundary conditions \eqref{soft} and $\nu\wedge\vE|_{\partial\Omega}=\nu\wedge\vE_0|_{\partial\Omega}$, we have $(\vE, \vH)$ satisfying
\begin{equation}\label{Iden_2}
-I=\int_{\Omega\setminus\bar{D}}\mu^{-1}|\nabla\wedge \vE-\nabla\wedge\vE_0|^2-
\omega^2\eps|\vE-\vE_0|^2dx+\int_D\mu^{-1}|\nabla\wedge\vE_0|^2-\omega^2\eps|\vE_0|^2dx
\end{equation}
where \[I\deq i\omega\int_{\partial\Omega}(\nu\wedge\vE_0)\cdot\left[\overline{(\nu\wedge \vH-\nu\wedge \vH_0)\wedge\nu}\right]dS.\]
\end{lemma}
Therefore, the proof essentially follows the "magnetic-hard" obstacle case, except that to show
\[\lim_{\tau\rightarrow\infty}\int_{\Osd}\eps|\vE-\vE_0|^2dx=0,\] we implement the regularity of the solution for an auxiliary boundary value problem similar to \eqref{auxbvp}
\[
\left\{\begin{array}{l}
\nabla\wedge\vP=i\omega\mu\vQ,\;\;\;\;\nabla\wedge\vQ=-i\omega\eps\vP-\frac{i}{\omega}\overline{\tvE}\;\;\;\;\mbox{ in }\;\Osd,\\
\nu\wedge\vP|_{\partial\Omega}=0,\\
\nu\wedge\vP|_{\partial D}=0,
\end{array}\right.
\]
see \cite{KSU}, and the Sobolev embedding theorem.

\vspace{0.5cm}

\begin{appendices}
\section{Well-posedness of a mixed boundary value problem for Maxwell's equations}
Let $\Omega\subset\R^3$ be a bounded domain with smooth boundary and $\overline{D}\subset \Omega$, Consider the boundary value problem of the Maxwell equation
\begin{equation}\label{Max_A1}
\left\{\begin{array}{l} \nabla\wedge\vE=i\omega\mu\vH\;\;\;\; \nabla\wedge\vH=-i\omega\gamma\vE\;\;\mbox{ in }\;\Osd,\\
\nu\wedge\vE|_{\partial\Omega}=f,\\
\nu\wedge(\mu^{-1}\nabla\wedge\vE)|_{\partial D}=i\omega\nu\wedge\vH|_{\partial D}=g.\end{array}\right.
\end{equation}
Here $\eps$ and $\gamma$ are complex-valued functions in $C^k(\Osd)$ with positive real parts and $\omega\in \C$.
\begin{theorem}\label{thm_A1}
There is a discrete subset $\Sigma$ of $\C$ such that for $\omega$ not in $\Sigma$, there exists a unique solution $(\vE, \vH)\in H^k_{\Div}(\Osd)
\times H^k_{\Div}(\Osd)$ of \eqref{Max_A1} given any $f\in TH^{k-1/2}_{\Div}(\partial\Omega)$ and $g\in TH^{k-1/2}_{\Div}(\partial D)$. The solution satisfies
\begin{equation}\label{A_5}\|\vE\|_{H^k_{\Div}(\Osd)}+\|\vH\|_{H^k_{\Div}(\Osd)}\leq C(\|f\|_{TH^{k-1/2}_{\Div}(\partial\Omega)}+\|g\|_{TH^{k-1/2}_{\Div}(\partial D)})\end{equation}
with $C>0$ independent of $f$ and $g$.
\end{theorem}

Here we proceed to prove the theorem by modifying the variational method in \cite{KSU}, \cite{C} and \cite{L}. From the Maxwell's equation \eqref{Max_A1}, the electric field $\vE$ satisfies the second
order equation
\[
\nabla\wedge(\mu^{-1}\nabla\wedge\vE)-\omega^2\gamma\vE=0\;\;\mbox{ in }\;\Osd,
\]
and
\begin{equation}\label{A_2}
\nabla\cdot\gamma\vE=0\;\;\mbox{ in }\;\Osd.
\end{equation}
Therefore, we consider
\begin{equation}\label{A_3}
\nabla\wedge(\mu^{-1}\nabla\wedge\vE)-s\gamma^{-1}\nabla(\nabla\cdot\gamma\vE)-\omega^2\gamma\vE=0\;\;\mbox{ in }\Osd
\end{equation}
where $s$ is a positive real number. The equation \eqref{A_2} will follow from
\[s\nabla\cdot(\gamma^{-1}\nabla(\nabla\cdot\gamma\vE))+\omega^2\nabla\cdot\gamma\vE=0\] which is obtained by taking divergence of
\eqref{A_3}.
Denoting the $L^2(\Osd)$ inner product by $(\cdot , \cdot)$ and $L^2(\Gamma)$ inner products by $\langle\cdot, \cdot\rangle_\Gamma$ (where $\Gamma=\partial\Omega$ or $\partial D$), we define the bilinear form associated with the elliptic system \eqref{A_3}:
\begin{equation}\label{A_bl}
B(\vE, \vF)\deq (\mu^{-1}\nabla\wedge\vE, \nabla\wedge\vE)+s(\nabla\cdot\gamma\vE, \nabla\cdot\gamma\vE)
\end{equation}
for $\vE, \vF\in X$ where
\[X=\{\vF\in (H^1(\Osd))^3\,|\,\nu\wedge\vF|_{\partial\Omega}=0,\,\nu\cdot\gamma\vF|_{\partial D}=0\}.\]

By Green's formulae, we have that $B$ is related to the differential operator
\[P=\nabla\wedge(\mu^{-1}\nabla\wedge)-s\gamma^{-1}\nabla(\nabla\cdot\gamma)\] by
\begin{equation}\label{bil}
B(\vE, \vF)=(P\vE, \vF)-\langle\nu\wedge(\mu^{-1}\nabla\wedge\vE),\vF\rangle_{\partial(\Osd)}+\langle s\nabla\cdot\gamma\vE, \nu\cdot\gamma\vF\rangle_{\partial(\Osd)}
\end{equation}
for $\vE, \vF\in (H^1(\Osd))^3$.
Then for $\tilde{f}\in X'$ the weak formulation of the mixed boundary value problem
\[
\left\{\begin{array}{l}P\vE=\tilde{f}\;\;\mbox{ in }\;\Osd, \\ \nu\wedge\vE|_{\partial\Omega}=0,\\
\nu\wedge(\mu^{-1}\nabla\wedge\vE)|_{\partial D}=0.
\end{array}\right.
\]
is: Find $\vE\in X$ such that
\[
B(\vE, \vF)=(\tilde{f}, \vF) \;\;\;\mbox{ for all }\vF\in X.
\]
By \eqref{bil}, this implies the natural boundary condition 
\begin{equation}\label{nbc}
\nu\wedge(\mu^{-1}\nabla\wedge\vE)|_{\partial D}=0,\;\;\;\;
\nabla\cdot\gamma\vE|_{\partial\Omega}=0.
\end{equation}\\

To show the theorem, one first has for the homogeneous boundary conditions,
\begin{proposition}\label{prop_2}
Suppose $\gamma$ and $\mu$ are complex functions in $C^1(\Osd)$ with positive real parts, and let $s$ be a positive real number. There is a discrete set $\Sigma_s\subset\C$ such that
if $\omega$ is outside this set, then for any $f\in X'$ there exists a unique solution $\vE\in X$ of
\begin{equation}\label{A_4}
\nabla\wedge(\mu^{-1}\nabla\wedge\vE)-s\gamma^{-1}\nabla(\nabla\cdot\gamma\vE)-\omega^2\gamma\vE=\tilde{f}
\end{equation}
satisfying
\[\|\vE\|_{H^1(\Osd)}\leq C\|\tilde{f}\|_{X'}.\]
\end{proposition}
\emph{Proof.} It's sufficient to show that $B$ is bounded and coercive on $X$. It's clear that $B$ is bounded,
\[|B(\vE, \vF)|\leq C\|\vE\|_{H^1(\Osd)}\|\vF\|_{H^1(\Osd)}.\]
To show the coercivity, first we have
\[|B(\vE, \vE)|\geq c\|\nabla\wedge\vE\|_{L^2(\Osd)}+s\|\nabla\cdot\vE\|_{L^2(\Osd)}-C\|\vE\|_{L^2(\Osd)}\|\vE\|_{H^1(\Osd)}.\]
It can be shown that there is a Poincar\'e inequality for 1-forms in $X$ similar to that in \cite{T}
\[\|\vE\|_{H^1(\Osd)}^2\leq C(\|\vE\|_{L^2(\Osd)}^2+\|\nabla\wedge\vE\|_{L^2(\Osd)}^2+\|\nabla\cdot\vE\|_{L^2(\Osd)}^2),\;\;\;\;\; \vE\in X,\] we have
\[B(\vE, \vE)\geq c\|\vE\|_{H^1(\Osd)}^2-C\|\vE\|_{L^2(\Osd)}^2.\]\qed\\

To show the higher order regularity of solutions, we define for $k\geq2$,
\[X^k=\{\vF\in (H^k(\Osd))^3\,|\,\nu\wedge\vF|_{\partial\Omega}=\nabla\cdot\gamma\vF|_{\partial\Omega}=\nu\cdot\gamma\vF|_{\partial D}=\nu\wedge(\mu^{-1}\nabla\wedge\vF)|_{\partial D}=0\}.\]
Then we have
\begin{proposition}\label{prop_3}
Let $\gamma$ and $\mu$ be functions in $C^k(\Osd)$, $k\geq2$, with positive real parts, and let $s>0$. Suppose $\omega\notin\Sigma_s$, then for any $\tilde{f}\in (H^{k-2}(\Osd))^3$ the equation \eqref{A_4} has a unique solution $\vE\in X^k$ and
\[\|\vE\|_{H^k(\Osd)}\leq C\|\tilde{f}\|_{H^{k-2}(\Osd)}.\]
\end{proposition}
This can be proved by the same techniques in Section 5.9 of \cite{T} for the Hodge Laplacian.\\

\emph{Proof of Theorem \ref{thm_A1}:} As in \cite{KSU}, we take $\Sigma$ to be the set $\Sigma_1$ in previous propositions, then $s$ can be chosen such that $\omega\notin\Sigma_s$ (for more details see \cite{KSU}).

To show uniqueness, suppose $(\vE, \vH)\in H^1_{\Div}(\Osd)\times H^1_{\Div}(\Osd)$ solves \eqref{Max_A1} with $f=g=0$. One has
\[\nabla\wedge(\mu^{-1}\nabla\wedge\vE)=\omega^2\gamma\vE,\;\;\;\;\; \nabla\cdot\gamma\vE=0.\]
It follows that $\vE\in X$ (by the natural boundary conditions) is a solution of \eqref{A_4} with $\tilde{f}=0$, which implies $\vE=\vH=0$ by Proposition \ref{prop_2}.

For existence, given $f\in TH^{k-1/2}_{\Div}(\partial\Omega)$ and $g\in TH^{k-1/2}_{\Div}(\partial D)\subset TH^{k-3/2}_{\Div}(\partial D)$, we can find $\vE_0\in H^k_{\Div}(\Osd)$ with
\[\nu\wedge\vE_0|_{\partial\Omega}=f,\;\;\;\;\nu\wedge(\mu^{-1}\nabla\wedge\vE_0)|_{\partial D}=g\] such that the extension is bounded, namely
\[\|\vE_0\|_{H^k(\Osd)}\leq C(\|f\|_{H^{k-1/2}(\partial\Omega)}+\|g\|_{H^{k-3/2}(\partial D)})\] 
Suppose $\tvE\in X^k$ is a solution, given by Proposition \ref{prop_3}, of \eqref{A_4} with
\[\tilde{f}=-\nabla\wedge(\mu^{-1}\nabla\wedge\vE_0)+s\gamma^{-1}\nabla(\nabla\cdot\gamma\vE_0)+\omega^2\gamma\vE_0\in (H^{k-2}(\Osd))^3.\]
Notice that \[\|\tilde{f}\|_{H^{k-2}(\Osd)}\leq C\|\vE_0\|_{H^k(\Osd)}.\]
Then $\vE=\vE_0+\tvE\in (H^k(\Osd))^3 $ satisfies
\begin{equation}
\nabla\wedge(\mu^{-1}\nabla\wedge\vE)-s\gamma^{-1}\nabla(\nabla\cdot\gamma\vE)-\omega^2\gamma\vE=0.
\end{equation}
This implies $\nabla\cdot\gamma\vE=0$ by an earlier argument and particular choice of $s$. If we define $\vH=\frac{1}{i\omega\mu}\nabla\wedge\vE\in (H^{k-1}(\Osd))^3$, then we have
$(\vE, \vH)$ the solution of Maxwell's equation. 

Applying the same argument to $\vH$, which satisfies a second order elliptic system by eliminating $\vE$ from the original Maxwell's equation. By uniqueness, one has $\vH\in (H^k(\Osd))^3$.

The fact $\vE\in H^k_{\Div}(\Osd)$ is obtained by
\[\Div(\nu\wedge\vE|_{\partial(\Osd)})=-\nu\cdot \nabla\wedge\vE|_{\partial(\Osd)}=-i\omega\mu\nu\cdot\vH|_{\partial(\Osd)}\in (H^{k-1/2}(\partial(\Osd)))^3.\]

Finally, the estimate \eqref{A_5} is derived from
\begin{eqnarray}\label{est_A1}\|\vE\|_{H^k(\Osd)}&\leq &\|\vE_0\|_{H^k(\Osd)}+\|\tvE\|_{H^k(\Osd)}\nonumber\\
&\leq& C(\|\vE_0\|_{H^k(\Osd)}+\|\tilde{f}\|_{H^{k-2}(\Osd)})\nonumber\\
&\leq& C\|\vE_0\|_{H^k(\Osd)}\nonumber\\
&\leq & C(\|f\|_{H^{k-1/2}(\partial\Omega)}+\|g\|_{H^{k-3/2}(\partial D)})\\
&\leq&C(\|f\|_{H^{k-1/2}(\partial\Omega)}+\|g\|_{H^{k-1/2}(\partial D)}).\nonumber\end{eqnarray} The same computation applies to $\vH$,
\begin{eqnarray*}\|\vH\|_{H^k(\Osd)}&\leq& C(\|f\|_{H^{k-3/2}(\partial\Omega)}+\|g\|_{H^{k-1/2}(\partial D)})\\
&\leq&C(\|f\|_{H^{k-1/2}(\partial\Omega)}+\|g\|_{H^{k-1/2}(\partial D)}).\end{eqnarray*}Then
\begin{eqnarray*}& &\|\Div(\nu\wedge\vE|_{\partial(\Osd)})\|_{H^{k-1/2}(\partial(\Osd))}\\
&\leq& C(\|\Div(f)\|_{H^{k-1/2}(\partial\Omega)}+\|\nu\cdot\vH\|_{H^{k-1/2}(\partial D)})\\
&\leq& C(\|\Div(f)\|_{H^{k-1/2}(\partial\Omega)}+\|\vH\|_{H^k(\Osd)})\\
&\leq& C(\|\Div(f)\|_{H^{k-1/2}(\partial\Omega)}+\|f\|_{H^{k-1/2}(\partial\Omega)}+\|g\|_{H^{k-1/2}(\partial D)})\\
&\leq& C(\|f\|_{TH^{k-1/2}_{\Div}(\partial\Omega)}+\|g\|_{TH^{k-1/2}_{\Div}(\partial D)}).\end{eqnarray*}
The same estimate can be obtained for $\vH$.
\qed

\end{appendices}


\begin{bibdiv}
\begin{biblist}

\bib{BU}{article}{
AUTHOR = {A. L. Bukhgeim and G. Uhlmann},
TITLE={Recovering a potential from partial Cauchy data},
JOURNAL={Comm. PDE},
Volume={27},
YEAR={2002},
PAGES ={653-668},
}


\bib{C}{article}{
Author={A. P. Calder\'{o}n},
TITLE={On an inverse boundary value problem},
JOURNAL={Seminar on Numerical Analysis and its Applications to Continuum Physics, Soc. Brasil. Mat., R\'{i}o de Janeiro},
YEAR={1980}
PAGES ={65-73},}

\bib{COS}{article}{
Author={P. Caro; P. Ola and M. Salo},
Title={Inverse Boundary Value Problem for Maxwell's Equations with Local Data},
Year={preprint},
}

\bib{C}{article}{
Author={M. Costabel},
Title={A Coercive Bilinear Form for Maxwell's equations},
Journal={J. Math. Anal. Appl.},
Volume={157},
Year={1991},
Pages={527-541}
}


\bib{I}{article}{
 AUTHOR = {M. Ikehata},
     TITLE = {How to draw a picture of an unknown inclusion from boundary measurements: Two mathematical inversion algorithms},
   JOURNAL ={J. Inverse Ill-Posed Problems},
    VOLUME = {7},
      YEAR = {1999},
     PAGES = {255-271},
}

\bib{I2}{article}{
Author={M. Ikehata},
Title={Reconstruction of the support function for inclusion from boundary measurements},
Journal={J. Inv. Ill-Posed Problems},
Volume={8},
Year={2000},
Pages={367-378},
}

\bib{IINSU}{article}{
Author={T. Ide, H. Isozaki, S. Nakata, S. Siltanen, G. Uhlmann}, Title={Probing for electrical inclusions with complex spherical waves},
Journal={Comm Appl. Math}, Volume={60}, year={2007}, number={10}, pages={1415-1442} }

\bib{KSU}{article}{
Author={C. Kenig; M. Salo; G. Uhlmann},
Title={Inverse problems for the anisotropic Maxwell equations},
Year={preprint},
}

\bib{L}{book}{
 AUTHOR = {R. Leis},
     TITLE = {Initial boundary value problems in mathematical physics},
   Publisher = {Teubner, Stuttgart/Willey, New York},
          YEAR = {1986},
}

\bib{M}{book}{
Author={W. Mclean},
Title={Strongly elliptic systems and boundary integral equations},
Publisher={Cambridge University Press},
Year={2000},}

\bib{NY}{article}{
 AUTHOR = {G. Nakamura; K. Yoshida},
     TITLE = {Identification of a non-convex obstacle for acoustical scattering},
   JOURNAL = {J. Inv. Ill-Posed Problems},
    VOLUME = {15},
      YEAR = {2007},
     PAGES = {1-14},
}

\bib{OPS1}{article}{
 AUTHOR = {P. Ola; L. P\"{a}iv\"{a}rinta; E. Somersalo},
     TITLE = {An inverse boundary value problem in electrodynamics},
   JOURNAL = {Duke Math. J.},
    VOLUME = {70},
      YEAR = {1993},
    NUMBER = {3},
     PAGES = {617-653},
}

\bib{OS}{article}{
 AUTHOR = {P. Ola; E. Somersalo},
     TITLE = {Electromagnetic inverse problems and generalized sommerfeld potentials},
   JOURNAL ={SIAM J. Appl. Math.},
    VOLUME = {56},
      YEAR = {1996},
    NUMBER = {4},
     PAGES = {1129-1145},
}


\bib{S}{article}{
 AUTHOR = {M. Salo},
     TITLE = {Semiclassical pseudodifferetial calculus and the reconstruction of a magnetic field},
   JOURNAL = {Comm. PDE},
    VOLUME = {31},
      YEAR = {2006},
    NUMBER = {11},
     PAGES = {1639-1666},
}

\bib{SIC}{article}{
Author={E. Somersalo; D. Issacson; M. Cheney},
Title={A linearized inverse boundary value problem for Maxwell's equations},
Journal={J. Compup. Appl. Math.},
volume={42},
Year={1992},
Pages={123-136},
}

\bib{SU}{article}{
Author={J. Sylvester; G. Uhlmann},
Title={Global uniqueness for an inverse boundary value problem},
Journal={Ann. Math.},
Volume={125},
Year={1987},
Pages={153-169},}

\bib{T}{book}{
 AUTHOR = {M. E. Taylor},
     TITLE = {Partial differential equations I: Basic Theory},
   Publisher = {Springer},
          YEAR = {1999},
}

\bib{U}{article}{
Author={G. Uhlmann},
Title={Commentary on Calder\'onÕs paper: on an Inverse Boundary
Value Problem},
Journal={Selecta (papers of Alberto P. Calder\'on), edited by A. Bellow; C. E. Kenig; P. Malliavin},}

\bib{UW}{article}{
Author={G. Uhlmann; J. N. Wang},
Title={Reconstructing discontinuities using complex geometric optics solutions},
Journal={SIAM J. Appl. Math.},
Volume={68},
number={4},
year={2008},
pages={1026-1044},}

\end{biblist}
\end{bibdiv}
\end{document}